\documentclass[12pt,draft]{amsart}
\usepackage{amssymb,amsmath,amsthm}
\usepackage[english]{babel}

\newcommand{\s}[1]{\mathbf{#1}}
\newcommand{\set}[1]{\{#1\}}

\newtheoremstyle{thm}
  {9pt}{9pt}{\itshape}{}{\bfseries}{}{.5em}{}
\theoremstyle{thm}
\newtheorem{thm}{Theorem}
\newtheorem{cor}[thm]{Corollary}
\newtheorem{lemma}[thm]{Lemma}
\newtheorem{prop}[thm]{Proposition}

\newtheoremstyle{defin}
  {9pt}{9pt}{}{}{\bfseries}{}{.5em}{}
\theoremstyle{defin}

\newtheoremstyle{exm}
  {9pt}{9pt}{}{}{\scshape}{}{.5em}{}
\theoremstyle{exm}
\newtheorem{exm}[thm]{Example}

\newtheoremstyle{prf}
  {}{}{}{}{\itshape}{:}{.5em}{}
\theoremstyle{prf}
\newtheorem*{prf}{Proof}

\DeclareMathOperator{\hght}{ht}

\title{On Goulden-Jackson's determinantal expression for the immanant}
\author{Matja\v z Konvalinka}
\date{\today}

\begin{document}

\begin{abstract}
 In 1992, Goulden and Jackson found a beautiful determinantal expression for the immanant of a matrix. This paper proves the same result combinatorially.
\end{abstract}

\maketitle
 
\section{Introduction}

Take a matrix $A = (a_{ij})_{i,j=1}^m$, its characteristic polynomial $\chi_A(t)=\det(A-tI)$ and its eigenvalues $\omega_1,\ldots,\omega_m$. Vieta's formulas tell us that the elementary symmetric functions
$$e_i(\omega_1,\ldots,\omega_m)=\sum_{j_1<\ldots<j_i} \omega_{j_1}\cdots \omega_{j_i}$$
are easily expressible in terms of $a_{ij}$:
$$e_0 t^m - e_1 t^{m-1} + \ldots + (-1)^m e_m =(t-\omega_1)\cdots(t-\omega_m) = $$
$$=\det(tI-A)=\sum_{i=0}^m (-1)^i t^{m-i}\sum_{J \in \binom{[m]}i} \det A_J,$$
with $A_J=(a_{ij})_{i,j \in J}$, i.e.
$$e_i(\omega_1,\ldots,\omega_m)=\sum_{J \in \binom{[m]}i} \det A_J.$$
The complete homogeneous symmetric functions
$$h_i(\omega_1,\ldots,\omega_m)=\sum_{j_1 \leq \ldots \leq j_i} \omega_{j_1}\cdots \omega_{j_i}$$
are also easy. We know that
$$\sum_i h_i t^i = \frac{1}{\sum_i (-1)^i e_i t^i}=\frac 1{t^m \sum_i (-1)^i e_i t^{i-m}}=\frac 1{t^m \det(t^{-1}I - A)}=\frac 1{\det(I-tA)}$$
and by MacMahon Master Theorem \cite[page 98]{macmahon}
\begin{equation} \label{eq1}
 h_i(\omega_1,\ldots,\omega_m)= \sum_{\lambda} a_{\lambda'_1\lambda_1}a_{\lambda'_2\lambda_2}\cdots a_{\lambda'_i\lambda_i},
\end{equation}
where:
\begin{itemize}
 \item $\lambda=\lambda_1\cdots \lambda_i$ runs over all sequences of $i$ letters from $[m]$, and
 \item $\lambda'=\lambda'_1\cdots \lambda'_i$ is the non-decreasing rearrangement of $\lambda$.
\end{itemize}

Goulden and Jackson \cite{goulden} proved the following.

\begin{thm} \label{intro2}
 Denote $h_i(\omega_1,\ldots,\omega_m)$ by $\Delta_i$, choose a partition $\lambda=(\lambda_1,\ldots,\lambda_p)$ of $n$, and write $\s a_\pi=a_{1\pi(1)}a_{2\pi(2)}\cdots a_{n\pi(n)}$ for a permutation $\pi \in S_n$. Then
 \begin{equation} \label {intro1}
  [\s a_\pi] \det(\Delta_{\lambda_i-i+j})_{p \times p} = \chi^\lambda(\pi),
 \end{equation}
 where $\chi^\lambda$ is the irreducible character of the symmetric group $S_n$ corresponding to $\lambda$.
\end{thm} 

Goulden and Jackson's result is stated in the (clearly equivalent) language of immanants. Theirs is one of the many papers in the early 1990's that brought about fascinating conjectures and results on immanants; see for example \cite{goulden2}, \cite{green}, \cite{ss}, \cite{stembridge} for details and further references.

\medskip

We will give two more proofs of this result. The first gives a recursion that specializes to Murnaghan-Nakayama's rule, and the second is a simple combinatorial proof of a statement equivalent to \eqref{intro1}.

\section{A recursive proof of Theorem \ref{intro2}}

Since $\Delta_{\lambda_i-i+j}$ is equal to $h_i(\omega_1,\ldots,\omega_m)$, $\det(\Delta_{\lambda_i-i+j})$ can be expressed as
$$s_\lambda(\omega_1,\ldots,\omega_m)$$
by the Jacobi-Trudi identity (see \cite[Theorem 7.16.1]{stanley}). Here $s_\lambda$ is the Schur symmetric function corresponding to the partition $\lambda$.

\medskip

Note first that all the terms $\s a = a_{1*}\cdots a_{1*} a_{2*}\cdots$ of
$$s_\lambda(\omega_1,\ldots,\omega_m)=\det (e_{\lambda_i'-i+j})$$
are \emph{balanced} in the sense that each $i$ appears as many times among $a_{i*}$ as among $a_{*i}$ (this also follows from \eqref{eq1} and the formula $s_\lambda=\det(h_{\lambda_i-i+j})$); here we are using the dual Jacobi-Trudi identity \cite[Corollary 7.16.2]{stanley}.

\medskip

Suppose we want to find the coefficient of $\s a = a_{1*}\cdots a_{1*} a_{2*}\cdots$. Assume that $C=\set{1,\ldots,k}$, $D=\set{k+1,\ldots,m}$ and that $\s a$ does not contain $a_{ij}$ for $i \in C,j \in D$ or $i \in D,j \in C$. The coefficient of $\s a$ does not change if we put all $a_{ij}$ that do not appear in $\s a$ to $0$; the matrix $A$ then has a block diagonal form $A_1 \oplus A_2$, and if $\xi_1,\ldots,\xi_k$ are the eigenvalues of $A_1$ and $\zeta_{k+1},\ldots,\zeta_m$ are the eigenvalues of $A_2$, the eigenvalues of $A$ are $(\omega_1,\ldots,\omega_m)=(\xi_1,\ldots,\xi_k,\zeta_{k+1},\ldots,\zeta_m)$. By definition,
$$s_\lambda(\omega_1,\ldots,\omega_m),$$
is the sum of $\s \omega^T = \omega_1^{\alpha_1(T)}\omega_2^{\alpha_2(T)}\cdots$ over all semistandard Young tableau (SSYT) $T$ of shape $\lambda$; here $\alpha_i(T)$ is the number of $i$'s in $T$. See \cite[\S 7.10]{stanley} for definitions and details. In every such $T$, the numbers $1,\ldots,k$ form a SSYT of some shape $\nu \subseteq \lambda$, and the numbers $k+1,\ldots,m$ form a SSYT of shape $\lambda/\nu$. Therefore

\begin{equation} \label{rec1}
 [\s a] s_\lambda(\omega_1,\ldots,\omega_m) = \sum_{\nu \subseteq \lambda} [\s a_1] s_\nu(\xi_1,\ldots,\xi_k) [\s a_2] s_{\lambda/\nu}(\zeta_{k+1},\ldots,\zeta_n).
\end{equation}

Since $s_{\lambda/\nu}$ is homogeneous of degree $|\lambda|-|\nu|$, we can restrict the sum to partitions $\nu \vdash k$.

\medskip

The second term in the product can be calculated explicitly by using the reverse Jacobi-Trudi identity (for example for $\s a_2$ ``small'') and the first term can be calculated recursively.

\begin{exm}
  Let us calculate the coefficient of $a_{11}a_{12}a_{21}a_{22}^2a_{34}a_{43}$. By \eqref{rec1}, we have to find the coefficient of $a_{34}a_{43}$ in $s_{\lambda/\nu}$ for $\nu=32,311,221$.\\
  We have
  $$s_{322/32}(\zeta_3,\zeta_4) = s_2(\zeta_3,\zeta_4)=h_2(\zeta_3,\zeta_4)=a_{33}^2 + a_{33}a_{44} + a_{34}a_{43} + a_{44}^2,$$
  $$s_{322/311}(\zeta_3,\zeta_4) = s_{11}(\zeta_3,\zeta_4) = e_2(\zeta_3,\zeta_4) = a_{33}a_{44} - a_{34}a_{43},$$
  $$s_{322/221}(\zeta_3,\zeta_4) = s_{21/1}(\zeta_3,\zeta_4) = e_1^2(\zeta_3,\zeta_4)=a_{33}^2 + 2a_{33}a_{44} + a_{44}^2$$
  and therefore
  $$[a_{11}a_{12}a_{21}a_{22}^2a_{34}a_{43}] s_{322}(\omega_1,\ldots,\omega_4) = [a_{11}a_{12}a_{21}a_{22}^2] \left(s_{32}(\zeta_1,\zeta_2) -  s_{311}(\zeta_1,\zeta_2)\right).$$
  Furthermore,
  $$s_{32}(\zeta_1,\zeta_2)=\begin{vmatrix} e_2 & 0 & 0 \\ e_1 & e_2 & 0 \\ 0 & e_0 & e_1\end{vmatrix}=e_2^2e_1=(a_{11}a_{22}-a_{12}a_{21})^2(a_{11}+a_{22})$$
  and
  $$s_{311}(\zeta_1,\zeta_2)=\begin{vmatrix} 0 & 0 & 0 \\ 1 & e_1 & 0 \\ 0 & 1 & e_1 \end{vmatrix}=0.$$
  Therefore
  $$[a_{11}a_{12}a_{21}a_{22}^2a_{34}a_{43}] s_{322}(\omega_1,\ldots,\omega_4) = -2.$$
\end{exm}

Assume that $\s a_2$ is of the form $a_{k+1,k+2}a_{k+2,k+3}\cdots a_{m,k+1}=:b_1\cdots b_l$. The corresponding matrix $A_2$ is
\begin{equation} \label{rec2}
 \begin{pmatrix} 0 & b_1 & 0 & \ldots & 0 \\ 0 & 0 & b_2 & \ldots & 0 \\ \vdots & \vdots & \vdots &\ddots & \vdots \\ 0 & 0 & 0 & \ldots & b_{l-1} \\ b_l & 0 & 0 & \ldots & 0 \end{pmatrix}
\end{equation}
and its characteristic polynomial is $(-1)^l (t^l - b_1b_2\cdots b_l)$. If the zeros are denoted $\eta_1,\ldots,\eta_l$, then $e_0(\eta_1,\ldots,\eta_l)=1$, $e_l(\eta_1,\ldots,\eta_l)=(-1)^{l-1}b_1\cdots b_l$, $e_i(\eta_1,\ldots,\eta_l)=0$ for $1 \leq i \leq l-1$.

\medskip

Recall that a border strip tableau is a connected skew shape with no $2 \times 2$ square, and that the height $\hght$ of a border strip tableau is defined to be one less than the number of rows. The following result is well known; we include the proof for the sake of completeness.

\begin{lemma}
 If $\lambda=(\lambda_1,\ldots,\lambda_p)$ is a partition and $\nu \subseteq \lambda$ with $|\lambda|-|\nu|=l$, then:
 \begin{enumerate}
  \item $\lambda/\nu$ is a border-strip tableau if and only if $\lambda'/\nu'$ is a border-strip tableau;
  \item $\lambda/\nu$ is a border-strip tableau if and only if $\lambda_1+p=l+1$;
  \item $\lambda/\nu$ is a border-strip tableau if and only if $\lambda_i=\nu_{i-1}+1$ for $2 \leq i \leq p$.
 \end{enumerate}
\end{lemma}
\begin{prf}
 (1) is obvious. (2) The squares of a border-strip tableau $\lambda/\mu$ form a NE-path from $(1,p)$ to $(\lambda_1,1)$. Each such path has $p-1+\lambda_1-1+1$ squares. Conversely, a partition from $(1,p)$ to $(\lambda_1,1)$ with $\lambda_1+p-1$ squares must be a NE-path and hence its squares form a border-strip tableau. (3) A $2 \times 2$ square in lines $i-1,i$ implies that $\lambda_i \geq \nu_{i-1}+2$. It is clear that if the squares of a tableau $\lambda/\mu$ form a NE-path, we have $\lambda_i=\nu_{i-1}+1$.\qed
\end{prf}

\begin{prop}
 If the zeros of $(-1)^l (t^l - b_1b_2\cdots b_l)$ are denoted by $\eta_1,\ldots,\eta_l$, we have $s_{\lambda/\nu}(\eta_1,\ldots,\eta_l)=0$ unless $\lambda/\nu$ is a border-strip tableau, and $s_{\lambda/\nu}(\eta_1,\ldots,\eta_l)=(-1)^{\hght(\lambda/\nu)}b_1b_2\cdots b_l$ for a border-strip tableau $\lambda/\nu$.
\end{prop}
\begin{prf}
 We can assume that $\lambda_i \neq \nu_i$ and that $\nu_p=0$. The indices of the entries of
 $$s_{\lambda/\nu}=\det(e_{\lambda'_i-\nu'_j-i+j})_{\lambda_1 \times \lambda_1}$$
 are strictly increasing in rows; the largest possible index is $\lambda_1+p-1=l$; and the indices of the diagonal elements are $e_{\lambda'_i-\nu'_i}$ with $1 \leq \lambda'_i-\nu'_i \leq l$ (and $\lambda'_i-\nu'_i=l$ only in the trivial case $\lambda_1=l$, $\nu=\emptyset$). By the lemma and the fact that $e_i=0$ for $1 \leq i \leq l-1$, the matrix is $0$ on and above the diagonal unless $\lambda/\nu$ is a border-strip tableau. If it is a border-strip tableau, the elements at $(i,j)$, $i \leq j < \lambda_1$, are $0$ (in particular, all the entries in the first row are $0$ except the last, which is $e_l$), and by the lemma, the subdiagonal elements are $1$ and the element $(1,\lambda_1)$ is $(-1)^{l-1} b_1\cdots b_l$. Therefore $s_{\lambda/\nu}(\eta_1,\ldots,\eta_l)=(-1)^{l-1+\lambda_1-1}b_1b_2\cdots b_l=(-1)^{p-1}b_1b_2\cdots b_l$, with $p-1=\hght(\lambda/\nu)$.\qed
\end{prf}

In other words, by \eqref{rec1}, $[\s a_\pi] s_\lambda(\omega_1,\ldots,\omega_m)$ satisfy the Murnaghan-Nakayama's rule, see \cite[Theorem 4.10.2]{sagan}. Together with the fact that Murnaghan-Naka\-yama's rule completely determines the irreducible characters $\chi^\lambda$, this shows that
$$[\s a_\pi] s_\lambda(\omega_1,\ldots,\omega_m)=\chi^\lambda(\pi).$$

Note that this also gives us the coefficient of $\s a_\pi = a_{1\pi(1)}\cdots a_{n\pi(n)}$ in $p_\lambda(\omega_1,\ldots,\omega_m)$: we know that
$$p_\lambda = \sum_{\mu} \chi^\mu(\lambda) s_\mu$$
and hence
$$[\s a_\pi] p_\lambda(\omega_1,\ldots,\omega_m)=\sum_{\mu} \chi^\mu(\lambda) \chi^\mu(\pi)$$
is (by the orthogonality of the table of characters) equal to 
$$z_\lambda=1^{j_1}j_1! 2^{j_2} j_2! \cdots$$
if the partition corresponding to $\pi$ is $\lambda=\langle 1^{j_1}2^{j_2}\cdots\rangle$, and $0$ otherwise; see e.g.\hspace{-0.07cm} \cite[Proposition 7.17.6]{stanley}.

\section{A proof of Theorem \ref{intro2} via scalar product}

Let us find the coefficient of $\s a_\pi$ in $e_\lambda(\omega_1,\ldots,\omega_m)=e_{\lambda_1}e_{\lambda_2}\cdots e_{\lambda_p}$. If we pick $a_{i\pi(i)}$ from $e_{\lambda_1}$, we must also pick $a_{\pi(i)j}$ from $e_{\lambda_1}$ because every term in $e_{\lambda_1}$ is balanced. But since $a_{\pi(i)\pi^2(i)}$ is the only term of $\s a_\pi$ with $\pi(i)$ as the first index, we must have $j=\pi^2(i)$. In other words, each of the cycles of $\pi$ must be chosen from one of the $e_{\lambda_i}$'s. We know that for $J=\set{j_1 < j_2 < \ldots < j_i}$ and $\tau$ a permutation of $j_1,\ldots,j_i$, $a_{j_1\tau(j_1)}\cdots a_{j_i\tau(j_i)}$ appears in $e_i$ with coefficient equal to the sign of $\tau$.

\medskip

This reasoning implies that the coefficient of $\s a_\pi = a_{1\pi(1)}\cdots a_{n\pi(n)}$ with $\pi$ of cycle type $\mu=(\mu_1,\ldots,\mu_q)$ in $e_\lambda$ is equal to $\varepsilon_\mu R_{\mu \lambda}$, where
\begin{itemize}
 \item $\varepsilon_\mu$ is equal to $(-1)^{j_2+j_4+\ldots}$ for $\mu=\langle 1^{j_1}2^{j_2}\cdots \rangle$, and
 \item $R_{\mu \lambda}$ is the number of ordered partitions $(B_1,\ldots,B_p)$ of the set $\set{1,\ldots,q}$ such that
$$\lambda_j = \sum_{i \in B_j} \mu_i$$
for $1 \leq j \leq p$.
\end{itemize}

\medskip

But we know that if $\langle \cdot,\cdot \rangle$ is the standard scalar product in the space of symmetric functions defined by $\langle h_\lambda,m_\mu \rangle = \delta_{\lambda\mu}$ and $\omega$ is the standard (scalar product preserving) involution given by $\omega(h_\lambda)=e_\lambda$, then $p_\mu = \sum_\nu R_{\mu \nu} m_\nu$ and $\omega(p_\mu)=\varepsilon_\mu p_\mu$ (see \cite[\S 7.4 -- \S 7.9]{stanley}) imply
$$\langle e_\lambda,p_\mu \rangle = \langle \omega(e_\lambda),\omega(p_\mu)\rangle = \varepsilon_\mu \langle h_\lambda,p_\mu \rangle = \varepsilon_\mu R_{\mu\lambda}.$$

Since $e_\lambda$ form a vector-space basis of the space of symmetric functions and since both the scalar product and the operator $[\s a_\pi]$ are linear, we have proved the following.

\begin{prop}
 For any symmetric function $f$ and for $\pi$ a permutation of cycle type $\mu$, we have
 $$[\s a_\pi] f(\omega_1,\ldots,\omega_m) = \langle f,p_\mu \rangle$$
 In particular,
 $$[\s a_\pi] s_\lambda(\omega_1,\ldots,\omega_m)=\langle s_\lambda,p_\mu \rangle = \chi^\lambda(\mu).\eqno \qed$$
\end{prop}

The proposition of course also implies that
$$[\s a_\pi] p_\lambda(\omega_1,\ldots,\omega_m) = \langle p_\lambda, p_\mu \rangle = z_\lambda \delta_{\lambda\mu},$$
which is what we proved at the end of the last section.

\section{An application to irreducible characters of $S_n$}

It is worthwhile to note the following determinantal description of irreducible characters of the symmetric group.

\begin{cor}
 Let $\lambda,\mu$ be partitions of $m$ and define $f_0,\ldots,f_m$ via the formula
 $$(t^{\mu_1}-u_1)(t^{\mu_2}-u_2)\cdots (t^{\mu_q}-u_q) = f_0 t^m - f_1 t^{m-1} + \ldots \pm f_m.$$
 Then
 $$\chi^\lambda(\mu) = [u_1\cdots u_q] \det(f_{\lambda_i'-i+j}).$$
\end{cor}

\begin{exm}
 Let $\lambda=(2,2,2,1)$ and $\mu=(3,2,2)$. Then
 $$(t^3-u_1)(t^2-u_2)(t^2-u_3) = t^7 -(u_2+u_3)t^5-u_1 t^4 + u_2 u_3 t^3+(u_1u_2+u_1u_3)t^2-u_1u_2u_3$$
 and so $f_0 = 1$, $f_1 = 0$, $f_2 = -u_2-u_3$, $f_3 = u_1$, $f_4 = u_2 u_3$, $f_5 = -u_1u_2-u_1u_3$, $f_6 = 0$, $f_7 = u_1u_2u_3$.
 Hence
 $$\chi^{2221}(322) = [u_1 u_2 u_3] \begin{vmatrix} f_4 & f_5 \\ f_2 & f_3 \end{vmatrix} = [u_1 u_2 u_3] (f_4f_3 - f_2f_5)=1-2=-1. \eqno \qed$$
\end{exm}

\begin{prf}
 Take the permutation $\pi = (1,2,\ldots,\mu_1)(\mu_1+1,\mu_1+2,\ldots,\mu_1+\mu_2)\cdots$ and form the block diagonal $A=A_1\oplus \ldots \oplus A_q$ with blocks of the form \eqref{rec2} corresponding to cycles of $\pi$. For $\mu=(3,2,2)$, the permutation is $(123)(45)(67)$ and the matrix is
 $$\begin{pmatrix} 0 & a_{12} & 0 & 0 & 0 & 0 & 0 \\
 		   0 & 0 & a_{23} & 0 & 0 & 0 & 0 \\
 		   a_{31} & 0 & 0 & 0 & 0 & 0 & 0 \\
 		   0 & 0 & 0 & 0 & a_{45} & 0 & 0 \\
 		   0 & 0 & 0 & a_{54} & 0 & 0 & 0 \\
 		   0 & 0 & 0 & 0 & 0 & 0 & a_{67} \\
 		   0 & 0 & 0 & 0 & 0 & a_{76} & 0 \end{pmatrix}.$$
 The characteristic polynomial of $A$ is 
 $$(t^{\mu_1}-u_1)(t^{\mu_2}-u_2)\cdots (t^{\mu_q}-u_q)$$
 for $u_1 = a_{12}a_{23}\cdots a_{\mu_11}$, $u_2 = a_{\mu_1+1,\mu_1+2}a_{\mu_1+2,\mu_1+3} \cdots,a_{\mu_1+\mu_2,\mu_1+1}$, etc. In other words, if $\omega_1,\ldots,\omega_m$ are the eigenvalues of $A$, $e_i(\omega_1,\ldots,\omega_m)=f_i$. But then $\det(f_{\lambda_i'-i+j})=s_\lambda(\omega_1,\ldots,\omega_m)$ and equation \eqref{intro1} implies
 $$[u_1\cdots u_q] \det(f_{\lambda_i'-i+j}) = [\s a_\pi] s_\lambda(\omega_1,\ldots,\omega_m) = \chi^\lambda(\mu). \eqno \qed$$
\end{prf}

\bigskip

{\bf Acknowledgements.} The author would like to thank Igor Pak, Alexander Postnikov and Richard Stanley for helpful discussions and suggestions.

\bigskip

\bigskip

\bigskip

{\sc \scriptsize Department of Mathematics, Massachusetts Institute of Technology, Cambridge, MA 02139\\
\tt{konvalinka@math.mit.edu}\\
\tt{http://www-math.mit.edu/\~{}konvalinka/}}

\end{document}